# A Bijection between Unbalanced Dyck Path and NE Lattice Path

Yannan Qian

University of Exeter, Penryn Campus, Penryn, TR10 9FE, UK

Nanjing University of Information, Science and Technology,

Ningliu Road 219, Nanjing, 211544, China

June 6, 2023

**Abstract:** Lattice paths are important tools on solving some combinatorial identities. This note gives a new bijection between unbalanced Dyck path (a path that never reaches the diagonal of the lattice) and NE (North and East only) lattice path from (0,0) to (n,n) by several partial reflections.

## Introduction

The combinatorial identity $\sum_{i=0}^{n}\binom{2i}{i}\binom{2n-2i}{n-i} = 2^{2n}$ can be proved easily by matching up the coefficients of generating function [1]. However, the direct combinatorial proofs are not so obvious. According to Stanley [1], the first combinatorial proof was given by G. Hajos, 1830s. Sved [2] gave a view of lattice path to solve this problem in 1984. The basic idea is separate the NE-lattice path of $2n$ steps (which has $2^{2n}$ permutations) into a NE lattice path from (0,0) to (i, i), which has $\binom{2i}{i}$ permutations, and a path of length $2n - 2i$ from (i, i) that would never reach the diagonal $y = x$. Such paths are called unbalanced Dyck paths in this note because they have different numbers of N-steps and E-steps. The difficulty is to prove the unbalanced Dyck path of length $2k$ has $\binom{2k}{k}$ permutations. A natural thought is that there are some bijections between unbalanced Dyck paths and NE lattice paths. Sved [2] gave a bijection by cutting and replacing the paths. This note gives another bijection by several partial reflections.

## Notations

For convenience, set the diagonal as a new $y$ axis. Then the N-step becomes (1,1) and the E-step becomes (1,-1). Call them upstep and downstep. Now a lattice path should end on the line $y = 0$. An unbalanced Dyck path will not touch the line $y = 0$ expect at (0,0).

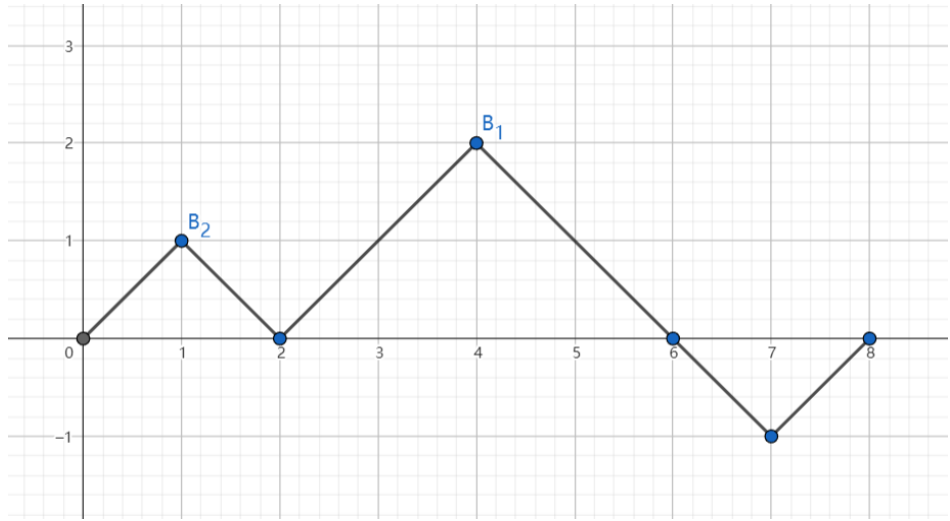

Figure 1

Suppose that the first step is an upstep. Note the highest point of the path $T$ as $B_1$. If there are several highest points, choose the left one. Then if there are any downsteps in front of $B_1$ (which means there are other peaks), note the highest point in front of $B_1$ as $B_2$, etc. Always choose the left one if there are points of the same height. We call a path down-Dyck path if it is a Dyck path that starts with a downstep. Note it as $D_K$ if it starts at point $K$. Use the notation $\overline{D_K}$ if the path is a reflection of $D_K$. Also, note down-unbalanced Dyck path starts at $K$ as $D'_K$.

Then a lattice path can be separated as
$$T = [upsteps, D_{B_n}, upsteps, D_{B_{n-1}}, \ldots, upsteps, D_{B_2}, upsteps, D'_{B_1}].$$

Figure 2 shows an example.

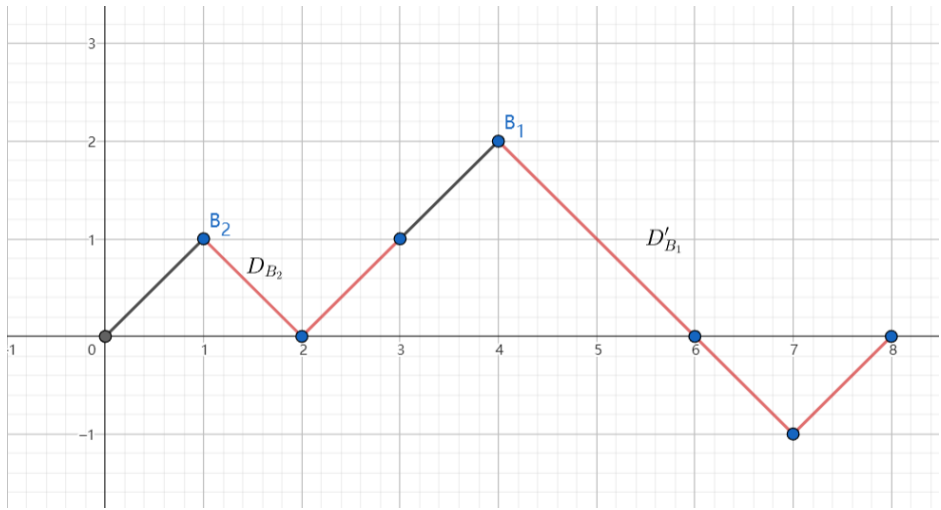

Figure 2

If the first step is not an upstep, we can do a reflection on $T$ and reflect it again after the bijection.

## The Bijection

Reflecting $D_{B_i}$ according to the line $y = y_{B_i}$, $D'_{B_1}$ according to the line $y = y_{B_1}$, and leaving upsteps in $T$ unchanged. Figure 3 shows an example.

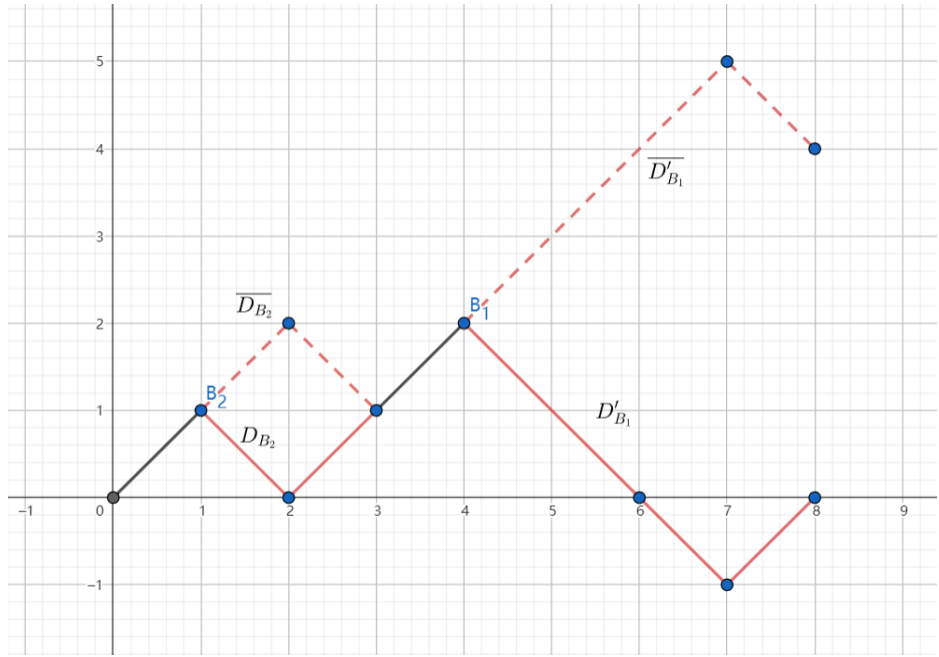

Figure 3

Then we have
$$\phi(T) = [upsteps, \overline{D_{B_n}}, upsteps, \overline{D_{B_{n-1}}}, \ldots, upsteps, \overline{D_{B_2}}, upsteps, \overline{D'_{B_1}}].$$
Obviously, $\phi(T)$ is an up-unbalanced Dyck path. Note the path $\phi(T)$ ends at point $G_1$. Then $y_{B_1} = \frac{y_{G_1}}{2}$.

Inverse:

1. Find the line $l: y = \frac{y_{G_1}}{2}$ to make sure $G_1$ will be on the line $y = 0$ after reflection.

2. Find the rightmost intersection $B_1$ of $l$ and $\phi(T)$. The intersection is the leftmost highest point of $T$. (All the points meet $\phi(T)$ but do not cross are not count, which indicate other highest points of $T$, but not the leftmost one.) Reflect path from $B_1$ to $G_1$ according to $l$.

3. Find the last downstep on the left of $B_1$. Note the end point of this downstep as $G_2$. $G_2$ should be the end point of the down-Dyck path $D_{B_2}$ in $T$, according to the separation of path $T$. Note the line $y = y_{C_2}$ as $l'$.

4. Find the rightmost intersection $l'$ and $\phi(T)$. That is $B_2$. Reflect path from $B_2$ to $G_2$ according to $l'$.

5. Do 3 and 4 by replacing the point $B$ and $G$ with the last $B$ and $G$ point found until there are no downsteps on the left of the last $B$.

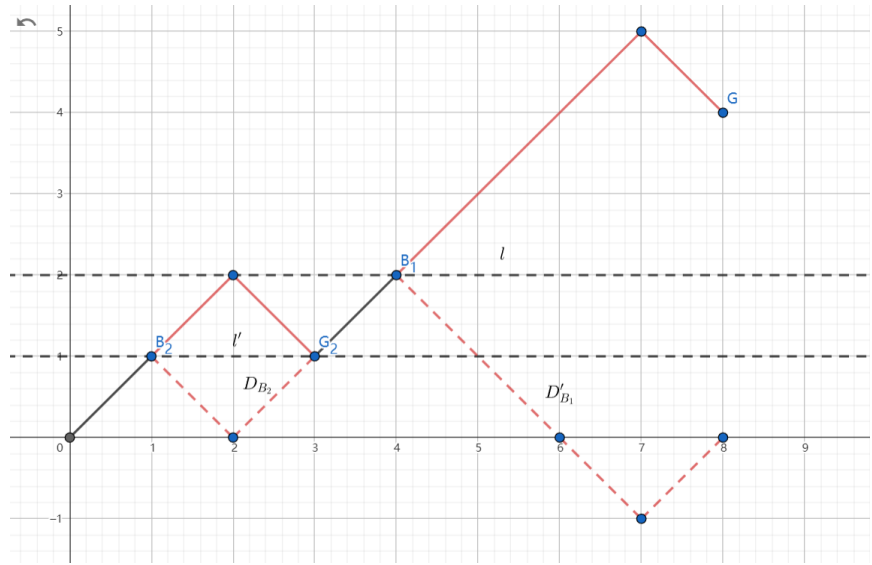

Figure 4

Both transforms are unique so it's a bijection. If the first step is a downstep, they are symmetry to the situations that the first step is an upstep.

If $T_1$ is not lower than $T_2$, then $\phi(T_1 + T_2) = \phi(T_1) + \overline{T}_2$, where the symbol '+' means to make the end point of the first part to be the start point of the second part.